\documentclass[12pt, leqno,twoside]{article}
\usepackage{amssymb}
\usepackage{amsmath}
\usepackage{amsthm}
\usepackage{graphicx}
\usepackage{float}
\usepackage{forest}
\usepackage{amsmath}
\usepackage{amsfonts}
\usepackage{amssymb}
\usepackage{amsmath,bbm,amssymb,amsxtra}
\usepackage{mathrsfs}
\usepackage{enumerate}
\usepackage{caption}
\allowdisplaybreaks
\usepackage{epsfig}
\usepackage{ae}

  \def\Xint#1{\mathchoice
    {\XXint\displaystyle\textstyle{#1}}%
    {\XXint\textstyle\scriptstyle{#1}}%
    {\XXint\scriptstyle\scriptscriptstyle{#1}}%
    {\XXint\scriptscriptstyle\scriptscriptstyle{#1}}%
    \!\int}

    \def\XXint#1#2#3{{\setbox0=\hbox{$#1{#2#3}{\int}$ }
    \vcenter{\hbox{$#2#3$ }}\kern-.6\wd0}}
    \def\dashint{\Xint-}

\newcommand{\vr}{{\varphi}}
\newcommand{\bx}{{\partial\Omega}}
\newcommand{\lbd}{{\lambda}}
\newcommand{\ms}{{\mathbb S}}
\newcommand{\md}{{\mathbb D}}

\newcommand{\diam}{\text{diam}}

\newcommand{\dist}{\text{\rm dist}}

\newcommand{\rarrow}{\rightarrow}
\newtheorem{thm}{Theorem}
\newtheorem{lem}{Lemma}
\newtheorem{prop}{Proposition}
\newtheorem{example}{Example}
\numberwithin{equation}{section}

\begin{document}
\title{\bf Controlled diffeomorphic extension of homeomorphisms
}
\author{Pekka Koskela, Zhuang Wang, and Haiqing Xu}
\date{ }
\maketitle

\begin{abstract}
Let $\Omega$ be an internal chord-arc Jordan domain and $\vr:\ms\rarrow\bx$ be a homeomorphism. We show that $\varphi$ has finite dyadic energy  if and only if $\vr$ has a diffeomorphic extension $h: \md\rarrow \Omega$ which has finite energy.

\end{abstract}

\section{Introduction}
Let $\Omega\subset\mathbb C$ be a bounded convex domain and suppose that $\varphi$ is a homeomorphism from the unit circle $\mathbb{S}$ onto $\partial \Omega.$ Then, by \cite{Kneser 1926 Jahresber. Deutsch. Math.-Verein.}, the 
complex-valued Poisson extension $h$ of $\varphi$ is a homeomorphism from $\bar{\mathbb{D}}$ onto $\bar{\Omega}$. This harmonic map $h$ is a diffeomorphism in $\md$ but its derivatives are not necessarily uniformly bounded. In $2007$, G. C. Verchota \cite{Verchota 2007 Proc. Amer. Math. Soc.} proved that the derivatives of $h$ may fail to be square integrable but that they are necessarily $p$-integrable over $\md$ for any $p<2$. In $2009$, T. Iwaniec, G. Martin and C. Sbordone improved on \cite{Iwaniec 2009 Discrete Contin. Dyn. Syst. Ser.} by showing that the derivatives belong to weak-$L^2$ with sharp estimates. In a related work \cite{Astala 2005 Proc. London Math. Soc.} by K. Astala, T. Iwaniec, G. Martin and J. Onninen, it was shown that if additionally $\bx$ is a $C^1$-regular Jordan curve, the square integrability of the derivatives of $h$ is equivalent to the requirement that
\begin{equation*}
\int_{\partial \Omega} \int_{\partial \Omega} |\log|\varphi^{-1}(\xi) -\varphi^{-1}(\eta)|| |d\xi| |d\eta| < +\infty.
\end{equation*}

In this note we give a generalization of the aforementioned results. Towards this end, recall that the Poisson extension of a homeomorphism $\vr: \ms\rarrow \partial \Omega$ may fail to be injective if $\Omega$ is not convex. Next, the boundary $\partial \Omega$ of a bounded convex domain $\Omega$ is a {\it chord-arc Jordan curve}: $\partial \Omega$ is a rectifiable Jordan curve and there is a constant $C$ such that for all $w_1, w_2\in \partial\Omega$,
$${\ell(w_1, w_2)}\leq C {|w_1- w_2|},$$
where $\ell(w_1, w_2)$ is the arc length of the shorter arc of $\partial\Omega$ joining $w_1$ to $w_2$. A domain whose boundary is a  chord-arc Jordan curve is called a {\it chord-arc Jordan domain}. Hence we are lead to ask for the optimal regularity of homeomorphic extensions $h: \md\rarrow \Omega$ for a given homeomorphism $\vr: \ms\rarrow \partial\Omega$ and a chord-arc Jordan domain $\Omega$.

Before introducing our first result, let us fix a dyadic decomposition $\{\Gamma_{j, k}: j\in \mathbb N, k=1, \cdots, 2^j\}$ of $\ms$, such that for a fixed $j\in \mathbb N$, $\{\Gamma_{j,k}: k=1, \cdots, 2^j\}$ is a family of arcs of length $2\pi/2^j$ with $\bigcup_{k}\Gamma_{j,k} = \ms$. The next generation is constructed in such a way that for each $k\in\{1, \cdots, 2^{j+1}\}$, there exists a unique number $k'\in \{1, \cdots, 2^j\}$, satisfying $\Gamma_{j+1,k}\subset \Gamma_{j,k'}$. Here, we call $\Gamma_{j, k'}$ the parent of $\Gamma_{j+1, k}$.  Fix $\Gamma_{1,1}$ to be the image of $[0, \pi]$ under the map $\theta\mapsto e^{i\theta}$. We denoted by $\ell(\vr(\Gamma_{j,k}))$  the arc length of the image arc of $\Gamma_{j, k}$ under the homeomorphism $\vr$.

\begin{thm}\label{thm1}
Let $\Omega\subset\mathbb C$ be an internal chord-arc Jordan domain and suppose that $\vr: \ms\rightarrow \partial\Omega$ is a homeomorphism. Then for $\lambda\in(-1, +\infty)$,  the following are equivalent:

(i) $\int_{\md}|Dh(z)|^2\log^\lambda(e+|Dh(z)|)\, dz<+\infty$ for some  diffeomorphic extension $h: \md\rarrow \Omega$ of $\vr$;

(ii) $\int_{\md} |Dh(z)|^2\log^\lambda(\frac{2}{1-|z|})\, dz<+\infty$ for some  diffeomorphic extension $h: \md\rarrow \Omega$  of $\vr$;

(iii) $\int_{\bx}\int_{\bx} |\log|\vr^{-1}(\eta)-\vr^{-1}(\xi)||^{\lambda+1}\, |d\eta|\, |d\xi|<+\infty;$

(iv) $\sum_{j=1}^{\infty}j^\lbd\sum_{k=1}^{2^j} \ell(\vr(\Gamma_{j,k}))^2<+\infty;$

(v) $\sum_{j=1}^{\infty}\sum_{k=1}^{2^j} \ell(\vr(\Gamma_{j,k}))^2\log^\lambda\left(e+\frac{\ell(\vr(\Gamma_{j, k}))}{2^{-j}}\right)<+\infty.$
\end{thm}


A Jordan domain $\Omega\subset\mathbb C$ is an  {\it internal chord-arc Jordan domain} if $\partial \Omega$ is rectifiable and there is a constant $C>0$ such that for all $w_1, w_2\in \partial\Omega$,
\begin{equation}\label{e1.1}
{\ell(w_1, w_2)}\leq C \lambda_\Omega(w_1, w_2),
\end{equation}
where $\ell(w_1, w_2)$ is the arc length of the shorter arc of $\partial\Omega$ joining $w_1$ to $w_2$ and $\lambda_\Omega(w_1, w_2)$ is the {\it internal distance} between $w_1, w_2$, which is defined as
 $${\lambda_\Omega(w_1, w_2)}=\inf_\alpha \ell(\alpha),$$
where the infimum is taken over all rectifiable arcs $\alpha \subset  \Omega$ joining $w_1$ and $w_2$; if there is no rectifiable curve joining $w_1$ and $w_2$, we let ${\lambda_\Omega(w_1, w_2)}=\infty$; cf. \cite[Section 3.1]{NV91} or \cite [Section 2]{OJB06}. 

Naturally, every chord-arc Jordan domain is also an internal chord-arc domain, but there are internal chord-arc domains that fail to be chord-arc; e.g., the inward cusp domain: 
$$\Omega_\epsilon = \md\setminus\{(x, y)\in \mathbb R^2: 0\leq x\leq 1, |y|\leq x^{2}\}.$$
 
 The statement of Theorem \ref{thm1} does not allow for $\lambda\leq -1$. Actually, (i)-(v) all hold for $\lambda<-1$, independently of $\vr$. When $\lambda=-1$, (iii) needs to be reformulated via a double logarithm, after which one still has a list of mutually equivalent conditions that may or may not hold, see \cite{KX}.
 
 In the case when $\Omega=\md$, Theorem \ref{thm1} can be formulated via the harmonic Poisson extension.

\begin{thm}\label{thm2}
Let $\vr: \ms\rarrow\ms$ be a homeomorphism and $h=P[\vr]: \md\rarrow\md$ be the harmonic Poisson extension of $\vr$. Then for $\lambda\in (-1, +\infty)$ the following are equivalent:

(i) $\int_{\md}|Dh(z)|^2\log^\lambda(e+|Dh(z)|)\, dz<+\infty;$

(ii) $\int_{\md} |Dh(z)|^2\log^\lambda(\frac{2}{1-|z|})\, dz<+\infty;$

(iii) $\int_{\ms}\int_{\ms} |\log|\vr^{-1}(\eta)-\vr^{-1}(\xi)||^{\lambda+1}\, |d\eta|\, |d\xi|<+\infty;$

(iv) $\sum_{j=1}^{\infty}j^\lbd\sum_{k=1}^{2^j} \ell(\vr(\Gamma_{j,k}))^2<+\infty;$

(v) $\sum_{j=1}^{\infty}\sum_{k=1}^{2^j} \ell(\vr(\Gamma_{j,k}))^2\log^\lambda\left(e+\frac{\ell(\vr(\Gamma_{j, k}))}{2^{-j}}\right)<+\infty.$

\end{thm}

Our main task is actually to prove Theorem \ref{thm2}. Indeed, once we know that Theorem \ref{thm2} holds, Theorem \ref{thm1} is obtained via a suitable change of variable, relying on the fact that there is a bi-Lipschitz map from $\bar \md$ onto $(\bar \Omega, \lambda_\Omega)$;  see Section \ref{s3.2}.

The paper is organized as follows. In Section 2 we give some relevant facts  about the dyadic decomposition of $\ms$ and some properties of the $N$-function $\Phi(t)= t^2\log^\lambda(e+t)$ for $\lambda>-1$.  Section 3 contains the full proofs of Theorem \ref{thm1} and Theorem \ref{thm2}. 

\section{Preliminaries}
\subsection{Dyadic decomposition}\label{s2.1}
Since our dyadic decomposition $\{\Gamma_{j, k}:  j\in \mathbb N, k=1, \cdots, 2^j\}$ of $\ms$ satisfies that $\Gamma_{1,1}$ is the image of $[0, \pi]$ under the map $\theta\mapsto e^{i\theta}$, we may assume that the dyadic decomposition $\{I_{j, k}=[2 \pi (k-1)/2^{j}, 2 \pi k/2^{j}]:  j\in \mathbb N, k=1, \cdots, 2^j\}$ of the interval $[0, 2\pi]$ matches with $\{\Gamma_{j, k}\}$ in the sense that each $\Gamma_{j, k}$ is exactly the image of $I_{j, k}$ under the map $\theta\mapsto e^{i\theta}$. Moreover, the dyadic arc $\Gamma_{j, k}$  is called a $j$-level dyadic arc, and we denote the two end points of $\Gamma_{j, k}$ by $\xi_{j,k}$ and $\xi_{j,k+1}$.   Two dyadic arcs which are $j$-level dyadic arcs for some $j\in \mathbb N$ are called brother dyadic arcs if they have a common parent.

\begin{figure}
 \centering
  \includegraphics[width=0.7\textwidth]{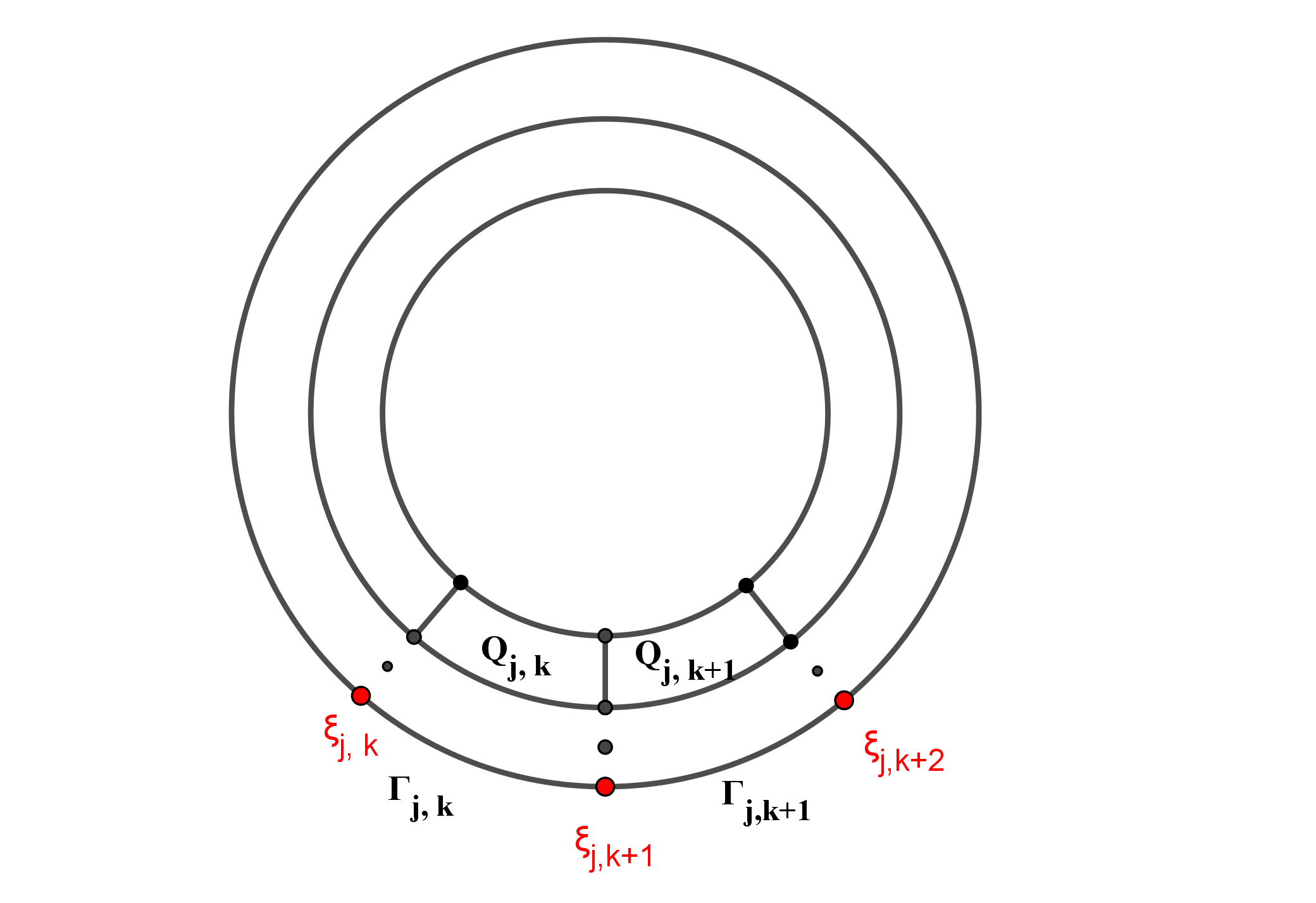}
   \caption{The decompositions of $\mathbb S$ and $\md$}
 \end{figure}

For the dyadic decompositions $\{\Gamma_{j, k}\}$ of $\ms$ and $\{I_{j, k}\}$ of $[0, 2\pi]$ as above, we have a decomposition of the unit disk $\md$ given by  $\{Q_{j,k}: j\in \mathbb N, k=1, \cdots, 2^j\}$ where
$$Q_{j, k}=\left\{re^{i\theta}: 1-1/2^{j-1}\leq r\leq 1-1/2^j, \theta\in I_{j, k}\right\}$$
for any $j\in\mathbb N$ and $k=1, \cdots, 2^j$, see Figure 1. It is easy to see that $\dist(Q_{j, k}, \ms)= 2^{-j}$ and that if $P: \md\rarrow \ms$ is the radial projection map, then $P(Q_{j, k})=\Gamma_{j, k}$ for all $j\in\mathbb N$ and $k=1, \cdots, 2^j$. 
There is an uniform constant $C>0$ such that for any $Q_{j,k}$ with center $x_{j,k}=r_{j, k}e^{i\theta_{j, k}}$, $r_{j, k}=1-3/2^{j+1}$ and $\theta_{j, k}=\pi(2k-1)/2^j$, we have 
\begin{equation*}
B(x_{j,k} , C^{-1} \diam(Q_{j,k})) \subset Q_{j,k} \subset B(x_{j,k} , C \diam(Q_{j,k})).
\end{equation*}
So for any $Q_{j,k},$ we can find a disk $B_{j,k}$ satisfying $B_{j,k} \subset Q_{j,k} \subset CB_{j,k},$ where $C$ is a constant independent of $Q_{j,k}.$

In order to later estimate the integral of $|Dh|$ over $Q_{j, k}$, we employ the following decomposition of $\ms$ via the dyadic arcs $\Gamma_{i, l}$ with $i\leq j$. The idea is to build a dyadic-type annular decomposition around $\Gamma_{j, k}$.

We fix $\Gamma_{j, k}$ and construct the decomposition via the following steps.
For simplicity, we assume that the unique brother arc of $\Gamma_{j,k}$ is $\Gamma_{j, k+1}$ and located at the anticlockwise side of $\Gamma_{j, k}.$

Step 1. On the anticlockwise side  of $\Gamma_{j,k+1},$ we choose the unique arc $\Gamma_{j-1, m_{j-1}^1}$ so that $\Gamma_{j, k+1}\cap \Gamma_{j-1, m_{j-1}^1}$ is a singleton. If the arc $\Gamma_{j-1, m_{j-1}^1 +1}$ shares the common parent with $\Gamma_{j-1, m_{j-1}^1}$, then set $(\Gamma_{j-1, m_{j-1}^1 +1})=\Gamma_{j-1, m_{j-1}^1 +1}$. If not, set $(\Gamma_{j-1, m_{j-1}^1 +1})=\emptyset$. Then define $\Gamma_{j-1, m_{j-1}^2}$ and $(\Gamma_{j-1, m_{j-1}^2 -1})$ analogously along the clockwise side of $\Gamma_{j, k}$.

Step 2. Repeat the process in Step 1 with $\Gamma_{j, k+1}$ replaced by $\Gamma_{j-1, m_{j-1}^1}$ if we have $(\Gamma_{j-1, m_{j-1}^1 +1})=\emptyset$ or $\Gamma_{j-1, m_{j-1}^1 +1}$ otherwise, and with $\Gamma_{j, k}$ replaced by $\Gamma_{j-1, m_{j-1}^2}$ if $(\Gamma_{j-1, m_{j-1}^2 -1})=\emptyset$ or $\Gamma_{j-1, m_{j-1}^2 -1}$ otherwise; unless 
we get $\Gamma_{j_0, m_{j_0}^1 }$, $(\Gamma_{j_0, m_{j_0}^1 +1})$, $\Gamma_{j_0, m_{j_0}^2 }$ and $(\Gamma_{j_0, m_{j_0}^2 -1})$ such that 
$$\left(\Gamma_{j_0, m_{j_0}^1 }\cup (\Gamma_{j_0, m_{j_0}^1 +1})\right) \bigcap \left(\Gamma_{j_0, m_{j_0}^2 }\cup (\Gamma_{j_0, m_{j_0}^2 -1})\right)\not=\emptyset.$$

This procedure decomposes $\ms$:
\begin{align*}
\mathcal P(\Gamma_{j,k})=\big\{  & \Gamma_{j, k}, \Gamma_{j, k+1}, 
\Gamma_{j-1, m_{j-1} ^{1}}, (\Gamma_{j-1, m_{j-1}^1 +1}), \Gamma_{j-1, m_{j-1}^2 }, (\Gamma_{j-1, m_{j-1}^2 -1}),  \\ 
& ..., \Gamma_{j_0, m_{j_0}^1 },(\Gamma_{j_0, m_{j_0}^1 +1}),\Gamma_{j_0, m_{j_0}^2 },(\Gamma_{j_0, m_{j_0}^2 -1}) \big\}. 
\end{align*}
where each $(\cdot)$ is either empty or a dyadic arc.
Correspondingly, we obtain a decomposition of $[0, 2\pi]$ as
\begin{align*}
\mathcal P(I_{j,k})=\big\{  & I_{j, k}, I_{j, k+1}, 
I_{j-1, m_{j-1} ^{1}}, (I_{j-1, m_{j-1}^1 +1}), I_{j-1, m_{j-1}^2 }, (I_{j-1, m_{j-1}^2 -1}),  \\ 
& ..., I_{j_0, m_{j_0}^1 },(I_{j_0, m_{j_0}^1 +1}), I_{j_0, m_{j_0}^2 },(I_{j_0, m_{j_0}^2 -1}) \big\}. 
\end{align*}

Let us illustrate this procedure by an example, see Figure 2 below. We use a binary tree to represent the dyadic decomposition $\{\Gamma_{j,k}\}$ of $\ms$. Suppose that we are given a dyadic arc $\Gamma_{4,7}$. Then $\Gamma_{4,8}$ is the unique brother arc of $\Gamma_{4,7}$. In Step 1, we choose $\Gamma_{3,5}$ and $\Gamma_{3.6}$ on the anticlockwise side of $\Gamma_{4,8},$ and $\Gamma_{3,3}$ on the clockwise side of $\Gamma_{4,7}.$ In Step 2, we choose $\Gamma_{2,4}$ on the anticlockwise side of $\Gamma_{3,6},$ and $\Gamma_{2,1}$ on the clockwise side of $\Gamma_{3,3}.$ There will be no more steps because $\Gamma_{2,1} \cap \Gamma_{2,4} \ne \emptyset.$ 

\begin{figure}[H]
    \centering \scalebox{0.72}
    {
    \begin{forest}
    [$\mathbb{S}^1$, for tree=
        [$\Gamma_{1,1}$,
            [$\boxed{\Gamma_{2,1}}$,
               [$\Gamma_{3,1}$,
                   [$\Gamma_{4,1}$]
                   [$\Gamma_{4,2}$]]
               [$\Gamma_{3,2}$,
                   [$\Gamma_{4,3}$]
                   [$\Gamma_{4,4}$]]]
            [$\Gamma_{2,2}$,
               [$\boxed{\Gamma_{3,3}}$,
                   [$\Gamma_{4,5}$]
                   [$\Gamma_{4,6}$]]
               [$\Gamma_{3,4}$,
                   [$\boxed{\Gamma_{4,7}}$]
                   [$\boxed{\Gamma_{4,8}}$]]]]
        [$\Gamma_{1,2}$,
            [$\Gamma_{2,3}$,
               [$\boxed{\Gamma_{3,5}}$,
                   [$\Gamma_{4,9}$]
                   [$\Gamma_{4,10}$]]
               [$\boxed{\Gamma_{3,6}}$,
                   [$\Gamma_{4,11}$]
                   [$\Gamma_{4,12}$]]]
            [$\boxed{\Gamma_{2,4}}$,
               [$\Gamma_{3,7}$,
                   [$\Gamma_{4,13}$]
                   [$\Gamma_{4,14}$]]
               [$\Gamma_{3,8}$,
                   [$\Gamma_{4,15}$]
                   [$\Gamma_{4,16}$]]]]
    ]
    \end{forest}
    }
    \caption{The decomposition of $\mathbb S$ around $\Gamma_{4,7}$}
\end{figure}
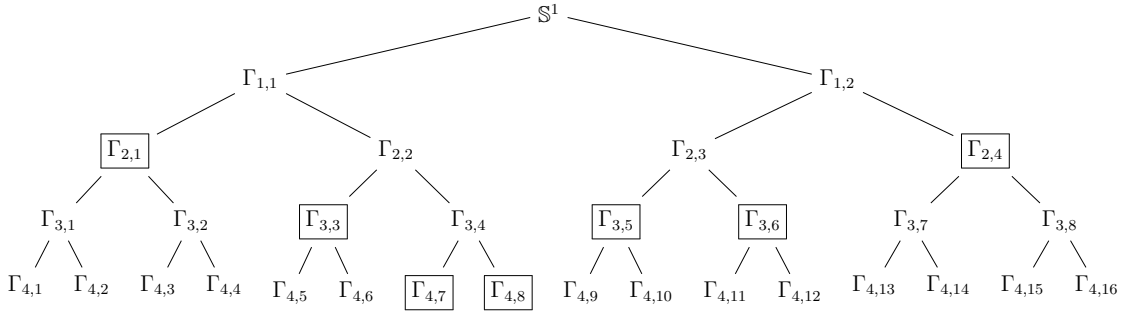

Suppose we are given an $n$-level arc $\Gamma_{n,m}.$ We are interested in the number of $j$-level arcs, $j\geq n$, which can induce $\Gamma_{n,m}$ by the above decomposition method.  Here, we say that a dyadic arc $\Gamma_{j, l}$  induces $\Gamma_{n, m}$ by the above decomposition method if and only if $\Gamma_{n, m}$ is in the set $\mathcal P(\Gamma_{j, l})$. 

To begin, $\Gamma_{n,m}$ can be first induced by its brother (without loss of generality assume that it is $\Gamma_{n,m+1}$), and two couples of $(n+1)$-level arcs. 
Notice that among all $(n+1)$-level arcs, there is only one couple that can induce both $\Gamma_{n,m}$ and $\Gamma_{n,m+1}.$ 
It follows that $\Gamma_{n,m}$ can be induced by three couples of $(n+1)$-level arcs.
Choose one of these $(n+1)$-level arcs. There is only one couple of $(n+2)$-level arcs which induce both $\Gamma_{n,m}$ and the chosen $(n+1)$-level arcs.
Then $\Gamma_{n,m}$ can be induced by six couples of $(n+2)$-level arcs.
Generally, for $j \geq n$ we have
\begin{equation}\label{number}
\sharp \{\Gamma: \Gamma \mbox{ is a }j\mbox{-level dyadic arc and  induce }\Gamma_{n,m}  \}\leq 3 \cdot 2^{j-n} .
\end{equation}

\subsection{$N$-functions}

A function $\Phi: [0, \infty)\rarrow [0, \infty)$ is an {\it $N$-function} if 
it is a continuous, increasing and convex function satisfying  $\Phi(0)=0$,
$$\lim_{t\rarrow 0+}\frac{\Phi(t)}{t}=0\ \ \text{and}\ \ \  \lim_{t\rarrow +\infty}\frac{\Phi(t)}{t}=+\infty.$$
An $N$-function $\Phi$ can be expressed as
$$\Phi(t)=\int_0^t \phi(s)\, ds,$$
where $\phi: [0, \infty)\rarrow [0, \infty)$ is an increasing, right-continuous function with $\phi(0)=0$ and $\underset{t \rightarrow + \infty}{\lim}\phi(t)=+\infty.$ 

For each $N$-function $\Phi$ and $t\geq 0$, set
$$\psi(t)=\sup_{\phi(s)\leq t} s\ \ \text{and}\ \ \Psi(t)=\int_0^t \psi(s)\, ds.$$
Then we call $\Psi$ the {\it complementary function} of $\Phi$. The complementary function of an $N$-function is also an $N$-function. We call $\Psi, \Phi$  a {\it pair of complementary $N$-functions}.

An $N$-function $\Phi$ is said to satisfy the $\Delta_2 -$condition if there is a constant $C_\Phi>0$, called a {\it doubling constant} of $\Phi$,  such that 
\begin{equation*}
\Phi(2 t) \le C_\Phi \Phi(t),\ \forall\ \  t \ge 0. 
\end{equation*}

\begin{prop}\label{doubling-p}
If an $N$-function $\Phi$ satisfies the $\Delta_2 -$condition, then for any constant $c>0$, there exist $c_1, c_2>0$ such that 
$$c_1\Phi(t)\leq \Phi(ct)\leq c_2\Phi(t)\ \ \ \ {\rm for\  all} \ \ \ t\geq 0,$$ 
where $c_1$ and $c_2$ depend only on $c$ and the doubling constant $C_\Phi$.  Therefore, we obtain that if $A\approx B$, then $  \Phi(A)\approx \Phi(B)$.
\end{prop}

\begin{lem}{\rm(\cite[Theorem 4.2]{KR61})}\label{doubling}
The complementary function of an $N-$function $\Phi$ satisfies the $\Delta_2$-condition on $[0,+\infty)$ if  there is a constants $l>1$ such that
\begin{equation*}
\Phi(t) \le \frac{1}{2l} \Phi(l t) \ \ \mbox{for any } t \ge 0. 
\end{equation*}  
\end{lem}

Given  an $N-$function $\Phi$, we denote by $L^{*} _{\Phi}$ the collection of all measurable functions $f$ such that $\int_{\mathbb{R}^n} \Phi(a f) <\infty$ for some $a >0.$

For a measurable function $f$ on $\mathbb R^2$, we define the Hardy-Littlewood maximal function of $f$  by setting
$$M_f(x)=\sup\dashint_B|f(z)|\,dz=\sup \frac{1}{|B|}\int_B |f(z)|\, dz,$$
where the supremum is taken over all open disks $B$ that contain $x$.

\begin{prop}{\rm{(\cite[Theorem 2.1]{Gallardo 1988 Publ. Mat.})}}\label{maximal-estimate}
Let $\Phi$ and $\Psi$ be a pair of complementary $N$-functions.
The following two conditions are equivalent:
\begin{enumerate}[(i)]
\item there exists positive constants $C$ and $b$ such that
\begin{equation*}
\int_{\mathbb{R}^n} \Phi(b M_f)(z) dz \le C \int_{\mathbb{R}^n} \Phi(|f|)(z) dz,\ \ \forall \ f \in L^{*} _{\Phi}, 
\end{equation*}
\item $\Psi$ satisfies the $\Delta_2$-condition on $[0,+\infty).$ 
\end{enumerate}
\end{prop}

\begin{example}\rm
Denote by $\Phi$ the function $\Phi(t) = t^2 \log^{\lambda} (e+t)$ for $\lambda>-1.$ An elementary computation shows that $\Phi$ is increasing, continuous and convex  on $[0,+\infty)$ with $\underset{t \rightarrow 0+}{lim} \frac{\Phi(t)}{t}=0$ and $\underset{t \rightarrow + \infty}{lim} \frac{\Phi(t)}{t}=+\infty.$ So $\Phi$ is an $N-$function. 
Moreover, both the function $\Phi$ and its complementary function satisfy the $\Delta_2$-condition. Hence we know that $\Phi$ satisfies Proposition \ref{doubling-p} and \ref{maximal-estimate}. We will use Proposition \ref{doubling-p} frequently in Section \ref{s3.1}.

Actually, a direct computation shows that the above $\Phi$ satisfies the $\Delta_2 -$condition on $[0,+\infty).$ 
 In order to check that the complementary $N-$function of $\Phi$ satisfies the $\Delta_2$-condition,  by Proposition \ref{doubling}, we only need to find a constant $l >1$ so that 
\begin{equation}\label{Phi+complementary delta2}
2 \log^{\lambda} (e+t) \le l \log^{\lambda}(e+lt),\ \forall \ t \ge 0.
\end{equation} 
In fact, if $\lambda \ge 0$ we can take $l =2.$ By monotonicity of $\log^{\lambda} (e+\cdot),$ we have that inequality \eqref{Phi+complementary delta2} holds for any $t \ge 0.$ 
If $-1< \lambda <0,$ we take $l=2^{1/(1+\lambda)}.$ 
By monotonicity we have
\begin{equation*}
\log (e+lt) \le l \log(e+t),\ \forall t \ge 0. 
\end{equation*}
Together with $l=\left(2/l\right)^{1/\lambda},$ it follows that \eqref{Phi+complementary delta2} holds for all $t \ge  0.$ 
\end{example}

\section{Proofs of Theorem \ref{thm1} and Theorem \ref{thm2}}
In this section, the notation $A\lesssim B$ means  that there is a constant  $C>0$ so that $A\leq C\cdot B$. Here and in this section, the notation $C$ denotes a positive constant  which may differ from line to line. The notation $A\approx B$ means  $A\lesssim B$ and $B\lesssim A$.

\subsection{Proof of Theorem \ref{thm2}}\label{s3.1}
\begin{proof}[{\bf Proof of (iv)$\Leftrightarrow$(v)}] 
Our claim is obvious when $\lambda=0$.

Case $\lbd>0$. Assume (iv) holds. Since $\ell(\vr(\Gamma_{j, k}))\leq 2\pi$ for any $\Gamma_{j, k}$ and $\lbd>0$, we have
$$\log^\lambda\left(e+\frac{\ell(\vr(\Gamma_{j, k}))}{2^{-j}}\right)\leq \log^\lambda\left(e+2^{j-1}\pi\right)\lesssim j^\lbd.$$
Therefore, we get
$$\sum_{j=1}^{\infty}\sum_{k=1}^{2^j} \ell(\vr(\Gamma_{j,k}))^2\log^\lambda\left(e+\frac{\ell(\vr(\Gamma_{j, k}))}{2^{-j}}\right)\lesssim\sum_{j=1}^{\infty}j^\lbd\sum_{k=1}^{2^j} \ell(\vr(\Gamma_{j,k}))^2<+\infty,$$
which gives the implication (iv)$\Rightarrow$(v). For the other direction, assume that (v) holds.  In order to estimate the  logarithmic term from below, we set
\begin{equation}\label{separate}
\chi(j, k)=\left\{\begin{array}{cc}
1, &\  \mathrm {if} \ \ \ell(\vr(\Gamma_{j, k})) >2^{-\frac{3j}{4}};\\
0,  & \mathrm{ otherwise}.
\end{array}
\right.
\end{equation}
Then
\begin{eqnarray*}
\sum_{j=1}^{\infty}j^\lbd\sum_{k=1}^{2^j} \ell(\vr(\Gamma_{j,k}))^2&=&\sum_{j=1}^{\infty}j^\lbd\sum_{k=1}^{2^j} \chi(j, k)\ell(\vr(\Gamma_{j,k}))^2\\
&&\ +\sum_{j=1}^{\infty}j^\lbd\sum_{k=1}^{2^j} (1-\chi(j, k))\ell(\vr(\Gamma_{j,k}))^2=: P_1+P_2.
\end{eqnarray*}
If $\ell(\vr(\Gamma_{j, k})) >2^{-\frac{3j}{4}}$, then we have
$$\log^\lambda\left(e+\frac{\ell(\vr(\Gamma_{j, k}))}{2^{-j}}\right)\geq \log^\lbd(e+2^{\frac 14 j})\gtrsim j^\lambda.$$
Hence, we obtain from (v) that
$$P_1\lesssim \sum_{j=1}^{\infty}\sum_{k=1}^{2^j} \ell(\vr(\Gamma_{j,k}))^2\log^\lambda\left(e+\frac{\ell(\vr(\Gamma_{j, k}))}{2^{-j}}\right)<+\infty.$$
For $P_2$, we have
$$P_2\leq \sum_{j=1}^{\infty}j^\lbd\sum_{k=1}^{2^j} 2^{-\frac32 j}=\sum_{j=1}^{+\infty}2^{-\frac12 j}j^\lambda<+\infty.$$
In conclusion, $P_1+P_2<+\infty$, and  (iv) follows. 

Case $\lbd<0$. Assume that (v) holds. Since $\ell(\vr(\Gamma_{j, k}))\leq 2\pi$ for any $\Gamma_{j, k}$ and $\lbd<0$, we have
$$\log^\lambda\left(e+\frac{\ell(\vr(\Gamma_{j, k}))}{2^{-j}}\right)\geq \log^\lambda\left(e+2^{j-1}\pi\right)\gtrsim j^\lbd.$$
Therefore, we get
$$\sum_{j=1}^{\infty}j^\lbd\sum_{k=1}^{2^j} \ell(\vr(\Gamma_{j,k}))^2\lesssim\sum_{j=1}^{\infty}\sum_{k=1}^{2^j} \ell(\vr(\Gamma_{j,k}))^2\log^\lambda\left(e+\frac{\ell(\vr(\Gamma_{j, k}))}{2^{-j}}\right)<+\infty,$$
which gives us the implication (v)$\Rightarrow$(iv). For the other direction, assume that (iv) holds.  Using the equation \eqref{separate}, if $\ell(\vr(\Gamma_{j, k})) >2^{-3 j/4}$, we have
$$\log^\lambda\left(e+\frac{\ell(\vr(\Gamma_{j, k}))}{2^{-j}}\right)\leq \log^\lbd(e+2^{\frac 14 j})\lesssim j^\lambda.$$
Therefore, we obtain
$$P_1':=\sum_{j=1}^{\infty}\sum_{k=1}^{2^j}\chi(j, k) \ell(\vr(\Gamma_{j,k}))^2\log^\lambda\left(e+\frac{\ell(\vr(\Gamma_{j, k}))}{2^{-j}}\right)\lesssim\sum_{j=1}^{\infty}j^\lbd\sum_{k=1}^{2^j} \ell(\vr(\Gamma_{j,k}))^2<+\infty.$$
 Moreover, since $\ell(\vr(\Gamma_{j, k})) \geq 0$ and $\lbd<0$, we always have
$$\log^\lambda\left(e+\frac{\ell(\vr(\Gamma_{j, k}))}{2^{-j}}\right)\leq 1.$$
Hence, we obtain
\begin{eqnarray*}
 P_2'&:=&\sum_{j=1}^{\infty}\sum_{k=1}^{2^j}(1-\chi(j, k)) \ell(\vr(\Gamma_{j,k}))^2\log^\lambda\left(e+\frac{\ell(\vr(\Gamma_{j, k}))}{2^{-j}}\right)\\
 &\leq& \sum_{j=1}^{\infty}\sum_{k=1}^{2^j} 2^{-\frac 32 j}= \sum_{j=1}^{\infty}2^{-\frac 12 j}<+\infty.
 \end{eqnarray*}
Then $P_1'+P_2'<+\infty$ which gives us (v). 

By combining the cases $\lambda=0$, $\lambda>0$ with $\lbd<0$ above, we finish the proof of  (iv)$\Leftrightarrow$(v).
\end{proof}

\begin{proof}[{\bf Proof of (iv)$\Rightarrow$(i)}]
The harmonic extension of $\vr$ is given by the Poisson integral formula:
$$h(z)=\frac{1}{2\pi}\int_{\ms}\frac{1-|z|^2}{|z-\zeta|^2}\vr(\zeta)\, |d\zeta|.$$
We can therefore compute the $z$-derivative by differentiating this kernel:
\begin{equation}\label{derivative}
h_z(z)=\frac{\partial h}{\partial z}(z)=\frac{1}{2\pi}\int_{\mathbb S}\frac{\zeta}{(z-\zeta)^2}\vr(\zeta)\, |d\zeta|=\frac{1}{2\pi}\int_0^{2\pi}\frac{e^{i\theta}}{(z-e^{i\theta})^2}\vr(e^{i\theta})\, d\theta.
\end{equation}

We may write
$$\vr(e^{i\theta})=\vr(1)\cdot e^{if(\theta)}$$
where $f: [0, 2\pi]\rarrow [0, 2\pi]$ is continuous and increasing with $f(0)=0$ and $f(2\pi)=2\pi$. This allows us to rewrite \eqref{derivative} as
\begin{equation}\label{Derivative}
h_z(z)=\frac{\vr(1)}{2\pi i}\int_0^{2\pi} \frac{\partial}{\partial \theta}\left[\frac{1}{(z-e^{i\theta})}\right]e^{if(\theta)}\, d\theta=\frac{\vr(1)}{2\pi i}\int_0^{2\pi}e^{if(\theta)} d\left(\frac{1}{(z-e^{i\theta})}\right).
\end{equation}
Here we use the Lebesgue-Stieltjes integration on the right-hand side of  \eqref{Derivative}. By integrating the right-hand side of  \eqref{Derivative} by parts, we obtain the estimate
$$h_z(z)=\frac{\vr(1)}{2\pi i}\int_0^{2\pi}\frac{1}{(z-e^{i\theta})}d\left(e^{if(\theta)}\right).$$
We therefore have 
$$|h_z(z)|\leq\frac{1}{2\pi}\int_0^{2\pi}\frac{1}{|z-e^{i\theta}|} \left|d\left(e^{if(\theta)}\right)\right|.$$
Denote by $\mu_f$ the Lebesgue-Stieltjes measure of the continuous increasing function $f$. Then we have $\left|d\left(e^{if(\theta)}\right)\right|\leq d(f(\theta))=d\mu_f(\theta)$. Therefore, we obtain the formula
\begin{equation}\label{est-der}
|h_z(z)|\leq\frac{1}{2\pi}\int_0^{2\pi}\frac{1}{|z-e^{i\theta}|}d\mu_f(\theta).
\end{equation}
Using our decomposition of the unit disk $\md$ and the decomposition $\mathcal P(I_{j,k})$ of $[0,2\pi]$ in Section \ref{s2.1}, for $\Phi(t)=t^2\log^\lbd(e+t)$, we get
\begin{align}\label{orlicz}
\int_{\md} \Phi(|h_z(z)|)\, dz= & \sum_{j=1}^{+\infty}\sum_{k=1}^{2^j}\int_{Q_{j, k}}\Phi(|h_z(z)|) dz \notag\\
\le & \sum_{j=1}^{+\infty}\sum_{k=1}^{2^j}\int_{Q_{j, k}} \Phi \left( \frac{1}{2 \pi} \sum_{I \in \mathcal P(I_{j,k})} \int_{I}  \frac{1}{|z-e^{i\theta}|}d\mu_f(\theta)\right) dz \notag\\
= & \sum_{j=1}^{+\infty}\sum_{k=1}^{2^j}\int_{Q_{j, k}} \Phi \left( \frac{1}{2 \pi} \sum_{n \le j} \sum_{m} \int_{I_{n,m}}  \frac{1}{|z-e^{i\theta}|}d\mu_f(\theta)\right) dz, 
\end{align}
where we abuse the notation: sum over those $m$ which belong to $\{1,\cdots, 2^n \}$ and satisfy $I_{n,m} \in \mathcal P(I_{j,k}).$ 
It is easy to check that $\#\{m\}\leq 3.$
For any $I_{n,m} \in \mathcal P(I_{j,k}),$ we know that $|z-e^{i\theta}|\approx 2^{-n}$ for any $z\in Q_{j, k}$ and $\theta\in I_{n, m}.$ Since $\mu_f$ is the Lebesgue-Stieltjes measure of our continuous increasing function $f$, for any interval $I_{n, m}$, we get $$\int_{I_{n, m}} d\mu_f(\theta) = |f(I_{n, m})|=\ell({\vr( \Gamma_{n, m})}).$$
Hence we obtain
\begin{align}\label{divide}
\int_{\md} \Phi(|h_z(z)|)\, dz \le & \sum_{j=1}^{+\infty}\sum_{k=1}^{2^j} \int_{Q_{j,k}}\Phi\left(\sum_{n\leq j}\sum_{m} \frac{\ell({\vr( \Gamma_{n, m})})}{2^{-n}}\right)dz \notag\\
\lesssim &\sum_{j=1}^{+\infty}\sum_{k=1}^{2^j}2^{-2j}\Phi\left(\sum_{n\leq j}\sum_{m} \frac{\ell({\vr( \Gamma_{n, m})})}{2^{-n}}\right) .
\end{align}
Since $\ell({\vr( \Gamma_{n, m})})\leq 2\pi$, we have
$$\sum_{n\leq j}\sum_{m} \frac{\ell({\vr( \Gamma_{n, m})})}{2^{-n}}\leq 2\pi \sum_{n\leq j}\sum_{m} \frac{1}{2^{-n}}\lesssim 2^j.$$
Using the same idea as in Proof of (iv)$\Leftrightarrow$(v), we have the estimate
$$\sum_{j=1}^{+\infty}\sum_{k=1}^{2^j}2^{-2j}\Phi\left(\sum_{n\leq j}\sum_{m} \frac{\ell({\vr( \Gamma_{n, m})})}{2^{-n}}\right)
\lesssim\sum_{j=1}^{+\infty}\sum_{k=1}^{2^j}2^{-2j}j^{\lambda} \left(\sum_{n\leq j}\sum_{m} \frac{\ell({\vr( \Gamma_{n, m})})}{2^{-n}}\right)^2+C.$$
This allows us to estimate \eqref{divide} as 
\begin{equation}\label{final-estimate}
\int_{\md} \Phi(|h_z(z)|)\, dz\lesssim \sum_{j=1}^{+\infty}\sum_{k=1}^{2^j}2^{-2j}j^{\lambda} \left(\sum_{n\leq j}\sum_{m} \frac{\ell({\vr( \Gamma_{n, m})})}{2^{-n}}\right)^2+C.
\end{equation}
For fixed ${\vr( \Gamma_{n, m})}$ and fixed $j\geq n$, the estimate \eqref{number} gives us $\#\{k\}\leq 3\cdot2^{j-n}.$
By applying the H\"older inequality on the right-hand side of inequality \eqref{final-estimate} and Fubini's theorem for series, we obtain
\begin{eqnarray*}
\int_{\md} \Phi(|h_z(z)|)\, dz
&\lesssim& \sum_{j=1}^{+\infty}\sum_{k=1}^{2^j}2^{-\frac 32j}j^{\lambda}\sum_{n\leq j}\sum_{m}\frac{\ell({\vr( \Gamma_{n, m})})^2}{2^{-\frac 32n}}+C\\
&\lesssim& \sum_{n=1}^{+\infty}\sum_{m=1}^{2^n} \frac{\ell({\vr( \Gamma_{n, m})})^2}{2^{-\frac 32n}}\left(\sum_{j\geq n}\sum_{k}2^{-\frac 32j}j^{\lambda} \right)+C\\
&\lesssim& \sum_{n=1}^{+\infty}\sum_{m=1}^{2^n} \frac{\ell({\vr( \Gamma_{n, m})})^2}{2^{-\frac 12n}}\begin{matrix} \underbrace{\left(\sum_{j\geq n}2^{-\frac 12j}j^{\lambda} \right)}\\ \approx 2^{-\frac 12 n}n^\lambda\end{matrix} +C\\
&\lesssim& \sum_{n=1}^{+\infty}\sum_{m=1}^{2^n}n^\lbd\ell({\vr( \Gamma_{n, m})})^2+C<+\infty.
\end{eqnarray*}

An analogous estimate follows for $h_{\bar z}$ by the same reasons. Thus we finish the proof of (iv)$\Rightarrow$(i). 
\end{proof}

\begin{proof}[{\bf Proof of (iv)$\Rightarrow$(ii)}]
Recall the estimate \eqref{est-der} for $|h_z|$ and formula  \eqref{orlicz}. From \eqref{orlicz} and the fact that $\log^\lambda(\frac{2}{1-|z|})\approx j^\lambda$ for any $z\in Q_{j, k}$, we have
\begin{eqnarray*}
\int_{\md}|h_z(z)|^2\log^\lambda\big(\frac{2}{1-|z|}\big)\, dz&=&\sum_{j=1}^{+\infty}\sum_{k=1}^{2^j}\int_{Q_{j, k}}|h_z(z)|^2\log^\lambda(\frac{2}{1-|z|})\, dz\\
&\approx&\sum_{j=1}^{+\infty}\sum_{k=1}^{2^j}j^\lambda\int_{Q_{j, k}}|h_z(z)|^2\, dz.
\end{eqnarray*}
By replacing $\Phi(t)=t^2\log^\lambda(e+t)$ by $t^2$ in \eqref{divide}, we obtain
$$\int_{Q_{j, k}}|h_z(z)|^2\, dz \lesssim 2^{-2j}\Big(\sum_{n\leq j}\sum_{m} \frac{\ell({\vr( \Gamma_{n, m})})}{2^{-n}}\Big)^2.$$
Hence we have
$$\int_{\md}|h_z(z)|^2\log^\lambda\big(\frac{2}{1-|z|}\big)\, dz\lesssim \sum_{j=1}^{+\infty}\sum_{k=1}^{2^j}2^{-2j}j^{\lambda} \left(\sum_{n\leq j}\sum_{m} \frac{\ell({\vr( \Gamma_{n, m})})}{2^{-n}}\right)^2.$$
Consequently, we may apply the very same arguments that we used above to estimate the right-hand side of \eqref{final-estimate} so as to arrive at
$$\int_{\md}|h_z(z)|^2\log^\lambda\big(\frac{2}{1-|z|}\big)\, dz\lesssim \sum_{n=1}^{+\infty}\sum_{m=1}^{2^n}n^\lbd\ell({\vr( \Gamma_{n, m})})^2<+\infty.$$

One may similarly deal with $h_{\bar z}$  to finish the proof of (iv)$\Rightarrow$(ii).
\end{proof}


\begin{proof}[{\bf Proof of (iv)$\Leftrightarrow$(iii)}]
We first divide the double integral on $\mathbb{S} $ into two parts:
\begin{align}\label{(iv) rarrow(iii): decom I+II}
 & \int_{\mathbb{S} } \int_{\mathbb{S} } \left|\log |\varphi^{-1}(\xi) - \varphi^{-1}(\eta)| \right|^{\lambda+1}  |d \eta| |d \xi| \notag \\  
 = & \int_{\mathbb{S} } \int_{\{\eta \in \mathbb{S}  : 1 \leq |\varphi^{-1}(\xi) - \varphi^{-1}(\eta)| \le 2\}} + \int_{\mathbb{S} } \int_{\{\eta \in \mathbb{S}  :  0 \le |\varphi^{-1}(\xi) - \varphi^{-1}(\eta)| <1\}} = I +II.
\end{align}
Since $ \lambda>-1$, we have $0\leq|\log t |^{\lambda+1}\leq (\log2)^{\lbd+1}$ for $1\le t\leq 2$. 
Therefore we obtain  $0 \le I\le C.$ So our  job is to estimate $II.$ 

Fubini's theorem changes $II$ to
\begin{align}\label{(iv) rarrow(iii):II+Fubini} 
II =& \int_{\mathbb{S} } \int_{\{\eta \in \mathbb{S}  : |\varphi^{-1}(\xi) - \varphi^{-1}(\eta)| < 1\}}\int_{0} ^{1} \Lambda(t) \chi_{\{|\varphi^{-1}(\xi) - \varphi^{-1}(\eta)| \leq t<1 \}} dt |d \eta| |d \xi|  \notag\\
= & \int_{\mathbb{S} } \int_{0} ^{1} \ell( \{\eta \in \mathbb{S} : |\varphi^{-1}(\xi) - \varphi^{-1}(\eta)| \le t\}) \Lambda(t) dt |d \xi|,
\end{align}
where $\Lambda(t) = (\lambda+1)  {\log^{\lambda} \left( \frac{1}{t} \right) }/{t}.$ Here we used the fact that 
 $$|\log t|^{\lambda+1}= \log^{\lambda+1} \left(\frac{1}{t} \right) =\int_t^1 \Lambda(t) dt, \ \text {for any } \ 0<t< 1.$$  
Then we have
\begin{eqnarray}\label{(iv) rarrow(iii): <}
II 
&=&  \sum_{j=0} ^{+\infty} \int_{\mathbb{S} } \int_{ 2^{-j-1} } ^{ 2^{-j}}\ell(\{\eta: |\varphi^{-1}(\xi) - \varphi^{-1}(\eta)|\le t \})    \Lambda(t) d \xi dt. \notag\\
&\le& C+\sum_{j=1} ^{+\infty} \int_{\mathbb{S} } \ell(\{\eta: |\varphi^{-1}(\xi) - \varphi^{-1}(\eta)|\le 2^{-j} \})  |d \xi| \int_{ 2^{-j-1} } ^{ 2^{-j}} \Lambda(t) dt.
\end{eqnarray}
If $\xi \in \varphi(\Gamma_{j,k})$ and $|h^{-1}(\xi) - h^{-1}(\eta)|\le  2^{-j}$, then the arc length of the shorter arc from $h^{-1}(\xi)$ to $h^{-1}(\eta)$ is at most $2^{-j}\pi$ and hence  it follows that $\eta \in \cup_{n=k-1} ^{k+1} \varphi(\Gamma_{j,n}).$
This means that for $j\geq 1$, we have
\begin{align}\label{(iv) rarrow(iii): <part 1}
\sum_{k=1} ^{2^j} &\int_{\varphi(\Gamma_{j,k})} \ell(\{\eta: |\varphi^{-1}(\xi) - \varphi^{-1}(\eta)|\le  2^{-j}\})  |d \xi|  \notag \\
& \le  \sum_{k=1} ^{2^j} \ell\big(\cup_{n=k-1} ^{k+1} \varphi(\Gamma_{j,n})\big) \ell(\varphi(\Gamma_{j,k}) 
\lesssim   \sum_{k=1} ^{2^j} \ell(\varphi(\Gamma_{j, k}))^2 . 
\end{align}
By the mean value theorem, we have
\begin{equation}\label{(iv) rarrow(iii): <part 2}
\int_{ 2^{-j-1}} ^{ 2^{-j}} \Lambda(t) dt = \log^{\lambda+1}(2^{j+1})- \log^{\lambda+1}(2^{j}) \approx  j^{\lambda}, \ \ j \ge 1.
\end{equation}
Combining (\ref{(iv) rarrow(iii): decom I+II}), (\ref{(iv) rarrow(iii):II+Fubini}), (\ref{(iv) rarrow(iii): <}), (\ref{(iv) rarrow(iii): <part 1}) and (\ref{(iv) rarrow(iii): <part 2}), we therefore deduce 
\begin{equation*}
\int_{\mathbb{S} } \int_{\mathbb{S} } | \log |\varphi^{-1}(\xi) - \varphi^{-1}(\eta)||^{\lambda+1}  |d \xi| |d \eta| 
\lesssim \sum_{j=1} ^{+\infty} \sum_{k=1} ^{2^j} \ell(\varphi(\Gamma_{j, k}))^2 j^\lambda
+C.
\end{equation*}

We are left to deal with the converse direction.
We divide the interval $[0,1]$ in (\ref{(iv) rarrow(iii):II+Fubini}) into intervals $[\pi /8, 1]$ and $[ 2^{-j+1}\pi, 2^{-j+2}\pi ]$ for $j \ge 5.$ 
Then we have 
\begin{equation}\label{(iv) rarrow(iii): >}
II \ge \sum_{j=5} ^{+\infty}  \int_{\mathbb{S} } \ell(\{\eta: |\varphi^{-1}(\xi) -\varphi^{-1}(\eta)| \le   2^{-j+1}\pi \}) |d \xi| \int_{2^{-j+1}\pi} ^{2^{-j+2}\pi} \Lambda(t) dt.
\end{equation}
Given any $j \ge 5$ and $1 \le k \le 2^j,$
  the inequality
\begin{equation*}
|\varphi^{-1}(\xi) -\varphi^{-1}(\eta) | \leq \ell(\Gamma_{j, k})\leq \pi 2^{1-j}
\end{equation*}
holds for all $ \eta, \xi \in \varphi(\Gamma_{j,k}).$ 
Thus we have
\begin{equation}\label{(iv) rarrow(iii): > part 1}
\sum_{k=1} ^{2^j} \int_{\varphi(\Gamma_{j,k})}  \ell(\{ \eta: |h^{-1}(\xi) -h^{-1}(\eta) | \le  \pi 2^{1-j} \})  |d \xi| 
\ge  \sum_{k=1} ^{2^j} \ell(\varphi(\Gamma_{j, k}))^2.
\end{equation}
Finally, the estimates (\ref{(iv) rarrow(iii): decom I+II}), \eqref{(iv) rarrow(iii): <part 2}, (\ref{(iv) rarrow(iii): >})   and (\ref{(iv) rarrow(iii): > part 1}) yield
\begin{equation*}
\int_{\mathbb{S} } \int_{\mathbb{S} } | \log|\varphi^{-1}(\xi) - \varphi^{-1}(\eta)||^{\lambda+1}   |d \xi| |d \eta| +C
\gtrsim \sum_{j=1} ^{+\infty} \sum_{k=1} ^{2^j} \ell(\varphi(\Gamma_{j, k}))^2 j^\lambda.
\end{equation*}
\end{proof}
\begin{proof}[{\bf  Proof of (ii)$\Rightarrow$(iv)}]
Since $\vr$ is uniformly continuous, there exists $j_0\geq 3$ such that $\ell(\vr( \Gamma_{j,k}))\leq \pi/3$ whenever $j\geq j_0$ and $k = 1, \cdots, 2^j$. Fix such $j, k$ and let $\xi\in \Gamma_{j, k-1}$, $\eta\in\Gamma_{j, k+1}$ ($\Gamma_{j, 0}=\Gamma_{j, 2^j}$, $\Gamma_{j, 2^j+1}=\Gamma_{j, 1}$). Since $h\in C(\bar \md)$ is a diffeomorphism, we have that 
\begin{equation*}
|h(\xi) -h(\eta)| \le \int_{\gamma} |Dh| \, ds.
\end{equation*}
for any rectifiable curve $\gamma\subset \md$ joining $\xi$ to $\eta$. Moreover, $|h(\xi_{j,k}) - h(\xi_{j,k+1})| \le |h(\xi) -h(\eta)|$ because $\varphi$ is homeomorphic and $\ell\big(\varphi(\Gamma_{j,k-1} \cup \Gamma_{j,k} \cup \Gamma_{j,k+1})\big) \le \pi$. Denote by $\hat\xi$ the midpoint of the arc $\Gamma_{j, k}$. Let $t_1>0$ be the smallest $t$ for which  $\partial B(\hat\xi, t) \cap \Gamma_{j,k+1} \neq \emptyset$ and let $t_2$ be correspondingly the largest such $t$. Write $\gamma_t=\partial B(\hat\xi, t) \cap\md$ for $t_1\leq t\leq t_2$. Then there exists an absolute constant $C$ so that 
$$\underset{t \in [t_1, t_2]}{\bigcup} \gamma_{t} \subset C Q_{j,k} \cap \mathbb{D}.$$
 Now 
$$|h(\xi_{j,k}) - h(\xi_{j,k+1})| \le \int_{\gamma_t} |Dh| \,ds$$
for $t\in [t_1, t_2]$ and by integrating with respect to $t$ we obtain
\begin{equation}\label{ii 2 iv: int t1 t2}
(t_2 -t_1) |h(\xi_{j,k}) - h(\xi_{j,k+1})| \le \int_{t_1}^{t_2}\int_{\gamma_t} |Dh| \,ds\, dt \leq \int_{CQ_{j,k} \cap \mathbb{D}} |Dh(z)|\, dz.
\end{equation}

Notice that $t_2 - t_1$ is uniformly comparable to $2^{-j}\approx \ell(\Gamma_{j,k})$ when $j \ge j_0$ and $k \in \{ 1, \cdots, 2^j \}.$
It follows from \eqref{ii 2 iv: int t1 t2} that 
\begin{equation}\label{ii 2 iv: h(xi jk) - h(xi jk+1)}
|h(\xi_{j,k}) - h(\xi_{j,k+1})| \lesssim \frac{1}{\ell(\Gamma_{j,k})} \int_{CQ_{j,k} \cap \mathbb{D}} |Dh(z)| \,dz.
\end{equation}
Recall that we have $\ell(\varphi( \Gamma_{j,k})) \le \pi/3$. Thus $\ell(\varphi(\Gamma_{j,k}))$ is also uniformly comparable to $|\varphi(\xi_{j,k}) - \varphi(\xi_{j,k+1})|.$ 
We can therefore conclude from \eqref{ii 2 iv: h(xi jk) - h(xi jk+1)} that for any $j \ge j_0$ and $1 \le k \le 2^j,$ the following inequality holds:
\begin{equation}\label{(ii) rarrow(iv): length < int_Q}
\ell(\vr(\Gamma_{j,k})) \lesssim \frac{1}{\ell (\Gamma_{j,k})} \int_{C Q_{j,k} \cap \mathbb{D}} |D h(z)| \, dz.
\end{equation} 

Fix $p \in (1,2)$ and $q \in (2,+\infty)$ with $1/p + 1/q =1.$ By applying the H$\ddot{\mbox{o}}$lder inequality to (\ref{(ii) rarrow(iv): length < int_Q}), we have 
\begin{equation*}
\ell(\vr(\Gamma_{j,k}))^{p} 
\lesssim \frac{1}{\ell(\Gamma_{j,k})^{p}} \int_{C Q_{j,k}} H(z)dz \left(\int_{C Q_{j,k} \cap \mathbb{D}} \log^{-\frac{\lambda q}{2} } \left(\frac{1}{1-|z|} \right)dz \right)^{p/q},
\end{equation*} 
where $H(z) = |D h (z)|^{p} \log^{{\lambda p}/{2}}\left( \frac{1}{1-|z|}\right)\chi_{\mathbb{D}} (z) .$
Moreover by changing to polar coordinates we have the estimate
\begin{align*}
\int_{C Q_{j,k} \cap \mathbb{D}} \log^{-\frac{\lambda q}{2}} \left(\frac{1}{1-|z|} \right)\, dz &
\lesssim \ell(\Gamma_{j,k}) \int_{0} ^{C \ell(\Gamma_{j,k})} \log^{-\frac{\lambda q}{2}} \left(\frac{1}{s} \right) ds \\
& \lesssim \ell(\Gamma_{j,k})^{2}\log^{-\frac{\lambda q}{2}} \left(\frac{1}{\ell(\Gamma_{j,k})} \right) \\
&\approx  \ell(\Gamma_{j,k})^{2} j^{-\lambda q/2}.
\end{align*}
It follows that
\begin{equation}\label{(ii) rarrow(iv):ell^p j^(lamda p/2) <}
\ell(\vr(\Gamma_{j,k}))^{p} j^{{\lambda p}/{2}} 
\lesssim \ell (\Gamma_{j,k})^{p} \dashint_{C Q_{j,k} \cap \mathbb{D}} H(z) dz.
\end{equation}
Next by the inclusion relationship between $B_{j,k}$ and $Q_{j,k}$ and the definition of Hardy-Litterwood maximal function we have
\begin{equation}\label{(ii) rarrow(iv): maximal}
\dashint_{C Q_{j,k}} H(z) dz 
\le  \dashint_{C B_{j,k}} H(z) dz  
\le \dashint_{B_{j,k}} M_{H} (z) dz 
\le \dashint_{Q_{j,k}} M_{H}  (z) dz.
\end{equation}
Combining (\ref{(ii) rarrow(iv):ell^p j^(lamda p/2) <}) with (\ref{(ii) rarrow(iv): maximal}) and then applying Jensen's inequality, we arrive at
\begin{equation*}
\ell(\vr(\Gamma_{j,k}))^{2} j^{\lambda} \lesssim \int_{Q_{j,k}} M_{H} ^{2/p}(z) dz, \ \  \forall j \ge j_0 \mbox{ and } 1 \le k \le 2^j.
\end{equation*} 
Then the $L^{2/p}$-boundedness of the Hardy-Littlewood maximal operator implies
\begin{align*}
\sum_{j=j_0} ^{+\infty} \sum_{k=1} ^{2^j} \ell(\vr(\Gamma_{j,k}))^{2} j^{\lambda} & 
\lesssim \sum_{j=1} ^{+\infty} \sum_{k=1} ^{2^j} \int_{Q_{j,k}} M_{H} ^{2/p} (z) \,dz \leq \int_{\mathbb R^2} M_{H} ^{2/p} (z) \,dz \\
& \lesssim \int_{\mathbb{D}} |D h(z)|^2 \log^{\lambda} \left(\frac{1}{1-|z|} \right) dz. 
\end{align*}
Moreover, we have that 
$$\sum_{j=1} ^{j_0 -1} \sum_{k=1} ^{2^j} \ell(\vr(\Gamma_{j,k}))^{2} j^{\lambda}\leq \sum_{j=1} ^{j_0 -1} j^{\lambda} \left(\sum_{k=1} ^{2^j} \ell(\vr(\Gamma_{j,k}))\right)^2 \lesssim  \sum_{j=1} ^{j_0 -1} j^{\lambda}\leq C.$$ 
 Hence we  conclude 
\begin{align*}
\sum_{j=0} ^{+\infty} \sum_{k=1} ^{2^j} \ell(\vr(\Gamma_{j,k}))^{2} j^{\lambda} = & \sum_{j=1} ^{j_0 -1} \sum_{k=1} ^{2^j} \ell(\vr(\Gamma_{j,k}))^{2} j^{\lambda} + \sum_{j=j_0} ^{+\infty} \sum_{k=1} ^{2^j} \ell(\vr(\Gamma_{j,k}))^{2} j^{\lambda} \\
\lesssim & C+ \int_{\mathbb{D}} |D h(z)|^2 \log^{\lambda} \left(\frac{1}{1-|z|} \right) dz. 
\end{align*}

\end{proof}



\begin{proof}[{\bf Proof of (i)$\Rightarrow$(v)}]   
Set $K(z)= |D h(z)| \chi_{\mathbb{D}}(z).$ As in the proof of (ii) $\Rightarrow$ (iv), we have 
\begin{equation*}
 \frac{\ell(\vr(\Gamma_{j,k}))}{\ell (\Gamma_{j,k})} \lesssim  \dashint_{C Q_{j,k} } K(z) dz\ \ \mbox{for any } j \geq j_0\mbox{ and } 1 \le k \le 2^j.
\end{equation*}
Thus Jensen's inequality for the function $\Phi(t) = t^2 \log^{\lambda} (e+t)$ and inequality (\ref{(ii) rarrow(iv): maximal}) with $H$ replaced by $K$ imply that
\begin{align}
\sum_{j=j_0} ^{+\infty} \sum_{k=1} ^{2^j} \ell(\Gamma_{j, k})^2 \Phi \left(\frac{\ell(\phi(\Gamma_{j,k}))}{\ell (\Gamma_{j,k})}\right)
&\lesssim  \sum_{j=1} ^{+\infty} \sum_{k=1} ^{2^j}\ell(\Gamma_{j, k})^2 \dashint_{ Q_{j,k} } \Phi(M_{K})(z) dz \notag\\
&\lesssim  \sum_{j=1} ^{+\infty} \sum_{k=1} ^{2^j}\int_{ Q_{j,k} } \Phi(M_{K})(z) dz \leq\int_{\mathbb{R}^2} \Phi(M_{K})(z) dz. \label{(v)<(i): < int_Q_j,k}
\end{align}
Here in the second inequality, we used the fact that $\ell(\Gamma_{j, k})^2\approx |Q_{j, k}|$.
Therefore Proposition \ref{maximal-estimate} gives 
\begin{equation}\label{(v)<(i):maximal+convex}
\int_{\mathbb{R}^2} \Phi(M_{K})(z) dz \lesssim \int_{\mathbb{R}^2} \Phi(K)(z) dz.
\end{equation}
Moreover, we have that 
$$\sum_{j=1} ^{j_0-1} \sum_{k=1} ^{2^j}  \ell(\vr(\Gamma_{j,k}))^2 \log^{\lambda}\left(e+\frac{\ell(\vr(\Gamma_{j,k}))}{\ell (\Gamma_{j,k})}\right) \lesssim\sum_{j=1} ^{j_0 -1} \sum_{k=1} ^{2^j} \ell(\vr(\Gamma_{j,k}))^{2} j^{\lambda}\lesssim  \sum_{j=1} ^{j_0 -1} j^{\lambda}\leq C.$$ 
Combining (\ref{(v)<(i): < int_Q_j,k}), (\ref{(v)<(i):maximal+convex}) with the above inequality gives
\begin{align*}
&\sum_{j=1} ^{+\infty} \sum_{k=1} ^{2^j}  \ell(\vr(\Gamma_{j,k}))^2 \log^{\lambda} \left(e+\frac{\ell(\vr(\Gamma_{j,k}))}{\ell (\Gamma_{j,k})}\right)  \lesssim C+\int_{\mathbb{D}} |Dh(z)|^2 \log^{\lambda}(e+|Dh(z)|) dz.
\end{align*}
\end{proof}

\subsection{Proof of Theorem \ref{thm1}}\label{s3.2}
\begin{proof}
From \cite[Lemma 3.2]{OJB06}, we know that any internal chord-arc Jordan domain $\Omega\subset \mathbb C$ is a bounded John disk whose boundary $\partial\Omega$ satisfies \eqref{e1.1}. Note also that in bounded John disks, the internal distance of any two boundary points is finite, \cite[Remark 6.6]{NV91}. 
By using the arc length parametrization of $\partial \Omega$ and the property \eqref{e1.1}, we see that there is a bi-Lipschitz map $g: \ms\rightarrow (\partial \Omega, \lambda_\Omega)$. Then applying \cite[Theorem 4.7]{OJB06} and the fact that the internal distance in the unit disk is the same as the Euclidean distance, we know that $g$ extends to a bi-Lipschitz map $\tilde g: \bar \md\rarrow (\bar\Omega, \lambda_\Omega)$. Moreover, the bi-Lipschitz map $\tilde g$ is a diffeomorphism in $\md$. Indeed, by the construction of the map $\tilde g$ in \cite[Theorem 4.7]{OJB06}, $\tilde g= \psi_1\circ \tilde f$, where $\psi_1: \md\rarrow \Omega$ is a conformal map and $\tilde f: \bar\md\rarrow \bar\md$ is a homeomorphic extension of $f=\psi_1^{-1}\circ g$ from \cite[Lemma 2.10]{GNV94}.
By the construction of $\tilde f$ in \cite[Lemma 2.10]{GNV94},
$\tilde f=\psi_2^{-1}\circ {\bar f}\circ \psi_2$, where $\psi_2: \bar \md\rarrow \bar {\mathbb H^2}$ is a M\"obius map and $\bar f: \bar {\mathbb H^2}\rarrow \bar {\mathbb H^2}$ is the Beurling-Ahlfors extension of the quasisymmetric map $\psi_2\circ f\circ \psi_2^{-1}$. Since conformal maps,  M\"obius maps and  Beurling-Ahlfors extensions of  quasisymmetric maps are diffeomorphic (cf. \cite[Chapter 8]{GMP17}), we obtain that the map $\tilde g$ is a diffeomorphism.

In the statements (i) and (ii) in Theorem \ref{thm1}, we let $h=\tilde g\circ P[{ g}^{-1}\circ \vr]$, where $P[{g}^{-1}\circ \vr]$ is the Poisson extension of ${g}^{-1}\circ \vr$. Hence $h$ is a diffeomorphism. Since $\tilde g$ is bi-Lipschitz with respect to the internal distance and the internal distance $\lambda_\Omega$ is the same as the Euclidean distance locally in $\Omega$,  we obtain that there is a constant  $C>0$ so that $1/C\leq |D{\tilde g}(z)|\leq C$ for any $z\in \md$. Hence the convergence of the integrals of $|Dh|$ in (i) and (ii) in Theorem \ref{thm1} is equivalent to  the statements (i) and (ii) for $P[{g}^{-1}\circ \vr]$ in Theorem \ref{thm2}. 

By \cite[Lemma 3.2, Lemma 3.4 and Lemma 3.6]{OJB06}, we know that for any $z_1, z_2\in \partial\Omega$, the shorter arc $\gamma_{z_1, z_2}$ of $\partial\Omega$ joining $z_1$ and $z_2$ satisfies
$$\diam(\gamma_{z_1, z_2}) \approx \lambda_{\Omega}(z_1, z_2).$$
Hence the arc length $\ell(\cdot)$ on $\partial \Omega$ with respect to the Euclidean distance is comparable to the arc length $\ell_{\lambda_\Omega}(\cdot)$ on $\partial \Omega$ with respect to the internal distance $\lambda_\Omega$. 
Combining with the bi-Lipschitz property of $g$, we obtain that $\ell(\vr(\Gamma_{j, k}))/C\leq \ell({g}^{-1}\circ\vr(\Gamma_{j, k}))\leq C\ell(\vr(\Gamma_{j, k}))$ for any $\Gamma_{j, k}$ and that we can use a change of variable via $g$ to conclude the equivalence of (iii) for $\vr^{-1}$ and $\vr^{-1}\circ g$ in Theorem \ref{thm1} and \ref{thm2}. Hence the statements (iii)-(v) for $\varphi$ in Theorem \ref{thm1} is equivalent to the statements (iii)-(v)  for ${g}^{-1}\circ \vr$ in Theorem \ref{thm2}. 
\end{proof}

\section*{Acknowledgment}
The first and second authors have been supported by the Academy of Finland via  Centre of Excellence in Analysis and Dynamics Research (project No. 307333).
The third author has been supported by the graduate school of University of Science and Technology of China, and by China Scholarship Council. This work was finished when the third author visited Department of Mathematics and Statistics at University of Jyv\"askyl\"a. He thanks the department for the hospitality. 

\bigskip
\bibliographystyle{amsplain}

\noindent Pekka Koskela

\noindent
Department of Mathematics and Statistics, University of Jyv\"askyl\"a, PO~Box~35, FI-40014 Jyv\"askyl\"a, Finland

\noindent{\it E-mail address}:  \texttt{pekka.j.koskela@jyu.fi}

\medskip
\noindent Zhuang Wang

\noindent
Department of Mathematics and Statistics, University of Jyv\"askyl\"a, PO~Box~35, FI-40014 Jyv\"askyl\"a, Finland

\noindent{\it E-mail address}:  \texttt{zhuang.z.wang@jyu.fi}

\medskip

\noindent Haiqing Xu

\noindent 
School of Mathematical Sciences, University of Science and Technology of China, Hefei 230026, P. R. China

\noindent{\it E-mail address}:  \texttt{hqxu@mail.ustc.edu.cn}

\end{document}